\documentclass[]{elsarticle}

\usepackage{hyperref}

\journal{Statistics \& Probability Letters}

\usepackage{numcompress}

\begin{document}

\newcommand{\nc}{\newcommand}
\nc{\nt}{\newtheorem}
\nt{defn}{Definition}
\nt{lem}{Lemma}
\nt{pr}{Proposition}
\nt{theorem}{Theorem}
\nt{cor}{Corollary}
\nt{ex}{Example}
\nt{ass}{Assumption}
\nt{step}{Step}
\nt{case}{Case}
\nt{subcase}{Subcase}
\nt{note}{Note}
\nc{\bd}{\begin{defn}} \nc{\ed}{\end{defn}}
\nc{\blem}{\begin{lem}} \nc{\elem}{\end{lem}}
\nc{\bpr}{\begin{pr}} \nc{\epr}{\end{pr}}
\nc{\bth}{\begin{theorem}} \nc{\eth}{\end{theorem}}
\nc{\bcor}{\begin{cor}} \nc{\ecor}{\end{cor}}
\nc{\bex}{\begin{ex}}  \nc{\eex}{\end{ex}}
\nc{\bass}{\begin{ass}}  \nc{\eass}{\end{ass}}
\nc{\bstep}{\begin{step}}  \nc{\estep}{\end{step}}
\nc{\bcase}{\begin{case}}  \nc{\ecase}{\end{case}}
\nc{\bsubcase}{\begin{subcase}}  \nc{\esubcase}{\end{subcase}}
\nc{\bnote}{\begin{note}}  \nc{\enote}{\end{note}}
\nc{\prf}{{\bf Proof.} }
\nc{\argmax}{\mathrm{argmax}}
\nc{\sgn}{\mathrm{sgn}}
\nc{\Var}{\mathrm{Var}}
\nc{\Cov}{\mathrm{Cov}}
\nc{\bak}{\!\!\!\!\!}

\newcommand\Binom[2]{\mbox{$\left(\!\!\!~^{~^{\textstyle
#1}}_{~_{\textstyle #2}}\right)$}}
\newcommand\binom[2]{\mbox{$\left(\!\!\!~^{~^{\scriptstyle
#1}}_{~_{\scriptstyle #2}}\right)$}}
\begin{frontmatter}

\title{Functional central limit theorems for the Nelson-Aalen and Kaplan-Meier estimators \\for dependent stationary data}

\author{Dragi Anevski}
\address{Centre for Mathematical Sciences, Lund University, Sweden}




\begin{abstract}
We derive process limit distribution results for the Nelson-Aalen estimator of a hasard function and for the Kaplan-Meier estimator of a distribution function, under different dependence assumptions. The data are assumed to be right censored observations of a stationary time series. We treat weakly dependent as well as long range dependent data, and allow for qualitative differences in the dependence for the censoring times versus the time of interest. \end{abstract}


\end{frontmatter}


\section{Introduction}
The estimation of a survival function based on censored survival data is of a fundamental importance in statistics. The Kaplan-Meier estimator is well known and well studied, and it's limit properties are well known, cf. $\!$Andersen et al. \cite{ABGK} for an introduction. The situation is similar to what is known for the empirical distribution function as an estimator of the distribution function, when we have independent data. However, that for the empirical distribution function, one has process limit distribution results also under general dependence assumptions, e.g. for mixing data and for long range dependent data, cf. $\!$Dehling and Taqqu \cite{DehlingTaqqu}. There are to our knowledge no process limit distribution results for the Kaplan-Meier estimator for dependent data. We establish such results in this paper.
\par
One approach to proving limit results for the Nelson Aalen estimator of the cumulative hasard function is to go via the theory of stochastic integration, expressing the Nelson-Aalen estimator as a stochastic integral of a previsible process with respect to a martingale, thus showing that in fact the Nelson-Aalen estimator is a martingale. Next the Rebolledo martingale central limit theorem is invoked to prove the weak convergence of the Nelson-Aalen estimator towards a Gaussian martingale, cf.$\!$ Anderson et al. $\!$\cite{ABGK}. The advantage of this approach is that it uses a very powerful theory to obtain strong results on e.g. $\!$variance estimates. The main disadvantage is that by relying on martingale results the approach is very restrictive in that it allows only independent data; weakly dependent data might potentially be managed with the theory of mixingales, cf.$\!$ McLeish \cite{McLeish}, but there is to our knowledge no functional central limit theorems for mixingales, cf. however Merlev\`ede and Peligrad \cite{MerlevedePeligrad} for functional central limit theorems for mixing data.
\par
We use the alternative approach to the derivation of the limit distribution of the Nelson-Aalen estimator which relies solely on two results: firstly a weak limit result of a certain bivariate empirical process, and secondly continuity of the map defining the Nelson-Aalen estimator, with respect to this bivariate process.  These two results will together, via the continuous mapping theorem, yield a shorter and more basic derivation of a limit result. There are several advantages over the martingale approach are: By relying on only empirical process theory we can apply the Nelson-Aalen estimator to a wider range of dependence situations, covering the known results for independent data but also implying new results for the marginal distributions for weakly dependent and long range dependent data. Furthermore the limit process can be written in an explicit form as compared to the usual reference to a Gaussian martingale with a certain mean and covariance function. A further advantage of this explicit representation is that it enables derivation of local and rescaled results; it is e.g.$\!$ easy to see what the process weak limit of the localized version $d_{n}^{-1}(\Lambda_{n}(t_{0}+sd_{n})-\Lambda_{n}(t_{0}))$ as $n\rightarrow\infty$ of the Nelson-Aalen estimator is, for $t_{0}$ fixed, $s$ varying and with $d_{n}\rightarrow 0$ as $n\rightarrow\infty$ a deterministic sequence. This may be useful when doing semiparametric estimation such as when using kernel smoothers or other smoothers for e.g. $\!$monotone hasard estimation, cf.$\!$ e.g.$\!$ Huang and Wellner \cite{HuangWellner}.
\par
The first limit process result for the Nelson-Aalen estimator was derived in Breslow and Crowley \cite{BreslowCrowley}, by using only empirical process theory together with smoothness of the functionals, i.e.$\!$ similar to our approach. A more recent paper with results similar to our results is Sun and Zhou \cite{SunZhou} that prove asymptotic normality of the Nelson-Aalen and Kaplan-Meier estimators, and for smoothed version of these, for weakly dependent ($\alpha$-mixing) life length data that are censored by i.i.d. data; however they state only pointwise limit distribution results.
\par
We will derive limit distribution results for the Kaplan-Meier estimate of the survival function, using the fact that the Kaplan-Meier estimator is a smooth map of the Nelson-Aalen estimator. It is in fact a compactly (or Hadamard) functionally differentiable map of the Nelson-Aalen estimator, cf.$\!$ Anderson et al.$\!$ \cite{ABGK}. Thus we can use the functional delta method, cf.$\!$ Gill \cite{GillHadamard}, to obtain the limit distributions of the Kaplan-Meier estimator. This generalises older results for independent and mixing data to also include long range dependent data.

The paper is organised as follows. In Section \ref{sec.2} we present the main results, that are functional central limit theorems for the Nelson-Aalen estimator (Theorem \ref{thm.1}), and for the Kaplan-Meier estimator (Theorem \ref{thm.2}).  In Section \ref{Sec:3} we specialise the results to qualitatively different dependence structures: Theorem \ref{thm.3} restates the well-known functional central limit theorem results for the i.i.d.$\!$ data case as well as state (new) results for stationary mixing data and Theorem \ref{thm.4} states (new) results on functional limit distributions for stationary long range dependent data. All proofs are gathered in  \ref{App}.

\section{Functional central limit theorems for the estimators; general dependence structures} \label{sec.2}
We are interested in making inference about the lifetime $T\geq 0$, from right-censored observations. We model the number of events that have occured by a counting process $N$, with intensity process $a(t)$ for which the multiplicative intensity model is assumed to hold, i.e.
\begin{eqnarray*}
   a(t)&=&\lambda(t) Y(t),
\end{eqnarray*}
where $\lambda(t)$ is the individual hasard and $Y(t)$ is the number at risk at time $t$. Then the cumulative hasard function
\begin{eqnarray*}
    \Lambda(t)&=&\int_{0}^{t}\lambda(u)\,du
\end{eqnarray*}
is of a central importance, for instance we have the relation $ \Lambda(t)=\log(1-F(t))$, where $F(t)=P(T\leq t)$ is the distribution function of the r.v. $T$.
\par
Assume the data are $(t_{i},\delta_{i}), i=1,\ldots,n$, where $t_{i}=\min(T_{i},\tau_{i})$ are the observed, possibly censored, life times, $\delta_{i}=1_{\{T_{i}\leq \tau_{i}\}}$ is an indicator for life or censoring event, $T_{i}$ are the partially unobserved lifetimes, while $\tau_{i}$ are the censoring times, and we assume we have independent and noninformative censoring, cf.$\!$ Anderson et al. $\!$\cite{ABGK}. Further assume that
\begin{eqnarray*}
  P(T_{i}\leq t)&=&F(t),\\
  P(\tau_{i}\leq t)&=&G(t),
\end{eqnarray*}
and define
\begin{eqnarray*}
     H(t)&=&P(\min(T_{i},\tau_{i})\leq t)=1-(1-F(t))(1-G(t)),\\
     H^{1}(t)&=&P(\min(T_{i},\tau_{i})\leq t,\delta_{i}=1),\\
     H^{0}(t)&=&P(\min(T_{i},\tau_{i})\leq t,\delta_{i}=0)=H(t)-H^{1}(t).
\end{eqnarray*}
Note that the (subdistribution) functions $H^{1}$ and $H^{0}$ are cadlag and increasing such that $H^{1}(t),H^{0}(t)\rightarrow 0$ as $t\rightarrow 0$. They are not distribution functions unless $P(\delta_{i}=1)=0$ or $1$, since $H^{0}(t)\rightarrow P(\delta_{i}=0)$ and $H^{1}(t)\rightarrow P(\delta_{i}=1)$ as $t\rightarrow\infty$, and never at the same time.
\par
Define also the empirical functions
\begin{eqnarray*}
   H_{n}(t)&=&\frac{1}{n}\sum_{i=1}^{n}1_{\{t_{i}\leq t\}},\\
   H_{n}^{1}(t)&=&\frac{1}{n}\sum_{i=1}^{n}\delta_{i}1_{\{t_{i}\leq t\}},\\
   H_{n}^{0}(t)&=&\frac{1}{n}\sum_{i=1}^{n}(1-\delta_{i})1_{\{t_{i}\leq t\}}=H_{n}(t)-H_{n}^{1}(t).
\end{eqnarray*}
Then the Nelson-Aalen estimator of the cumulative hasard function can be written as
\begin{eqnarray}
   \Lambda_{n}(t)&=&\int_{0}^{t}\frac{1}{1-H_{n}(u-)}\,dH_{n}^{1}(u),\label{eq:NA}
\end{eqnarray}
cf.$\!$ Shorack and Wellner \cite{Shorackwellner} and Anderson et al. $\!$\cite{ABGK}, and furthermore the cumulative hasard function can be written $ \Lambda(t)=\int_{0}^{t}\frac{1}{1-H(u-)}\,dH^{1}(u)$.
\par
In the case of independent censored data $\{t_{i}\}$, the Nelson-Aalen estimator can be interpreted as the nonparametric maximum likelihood estimator of $\Lambda$, i.e. $\! \Lambda_{n}$ maximizes
\begin{eqnarray}
      dP&=&\prod_{t} (Y(t)d\Lambda(t))^{dN(t)}(1-d\Lambda(t))^{Y(t)-dN(t)},\label{eq:likelihood}
\end{eqnarray}
cf.$\!$ Anderson et al.$\!$ \cite{ABGK}. If instead $\{t_{i}\}$ is a stationary sequence of dependent data, $(\ref{eq:likelihood})$ is not the full likelihood. However $(\ref{eq:likelihood})$ can still be interpreted as a partial likelihood. Furthermore, the function $\hat{\Lambda}$ that maximizes $(\ref{eq:likelihood})$ is obtained as an optimum over the same class of hasard functions as in the independent data case irrespective of what sort of dependence the data exhibits. This means that the Nelson-Aalen estimator can be seen as the (for dependent data ``marginal'') nonparametric maximum likelihood estimator of the survival function, no matter what dependence structure we have.
\par
In order to give a simple derivation of the limit distributions for the Nelson-Aalen estimator $(\ref{eq:NA})$, we will use the weak limits of the bivariate process $(H_{n},H_{n}^{1})$, properly centered and normalized, together with continuity results for the map
\begin{eqnarray}
     W: {\mathbf D}([0,\infty),{\mathbf R}) \times {\mathbf D}([0,\infty),{\mathbf R})  \ni (x,y)\mapsto \int_{0}^{\cdot}\frac{1}{1-x(u-)}\,dy(u)\in {\mathbf C}[0,\infty).\label{eq:W_map}
\end{eqnarray}
Let ${\mathbf D}[0,\infty)$ denotes the space of right continuous functions with left hand limits (cadlag) defined on $[0,\infty)$, equipped with the supnorm metric on compact intervals. Further we use the $\sigma$-algebra on ${\mathbf D}[0,\infty)$ generated by the closed balls in order to avoid measurability problems for the empirical process, cf.$\!$ Pollard \cite{Pollard}. Then ${\mathbf D}([0,\infty),{\mathbf R}) \times {\mathbf D}([0,\infty),{\mathbf R})$ is the set of bivariate cadlag functions, with metric equal to the maximum of the two supnorm on compact intervals metrics on ${\mathbf D}[0,\infty)$, so that a bivariate sequence $x_{n}=(x_{n}^{1},x_{n}^{2})$ converges iff on all compact intervals both $x_{n}^{1}$ and $x_{n}^{2}$ converge uniformly. 
\par
Define the bivariate stochastic process $w_{n}=(w_{n},w_{n}^{1})$, where
\begin{eqnarray*}
        w_{n}&=&a_{n}^{-1}(H_{n}-H),\\
        w_{n}^{1}&=&a_{n}^{-1}(H_{n}^{1}-H^{1}),
\end{eqnarray*}
where $a_{n}\downarrow 0$ as $n\rightarrow\infty$ is a deterministic sequence. Anticipating on the main results we make the following assumption.
\bass Assume that there are stochastic processes $w,w^{1}$ in ${\mathbf D}[0,\infty)$, such that
\begin{eqnarray*}
       (w_{n},w_{n}^{1})&\stackrel{\cal L}{\rightarrow}& (w,w^{1}),
\end{eqnarray*}
on ${\mathbf D}([0,\infty),{\mathbf R}) \times {\mathbf D}([0,\infty),{\mathbf R})$, as
$n\rightarrow\infty$.
 \eass 
We will derive the statement of Assumption 1 in the different applications by proving convergence of the finite-dimensional distributions of $w_{n}$, and tightness of $\{w_{n}\}_{n\geq 1}$. For the tightness the next result is useful.
\blem
Assume $\{w_{n}\}_{n\geq 1}$ and $\{w_{n}^{1}\}_{n\geq 2}$ are two sequences of processes that are separately tight on $\mathbf{D}[0,\infty)$. Then $\{(w_{n},w_{n}^{1})\}_{n\geq 1}$ is a sequence tight in ${\mathbf D}([0,\infty),{\mathbf R}) \times {\mathbf D}([0,\infty),{\mathbf R})$.
\elem

The form of $w,w^{1}$ and $a_{n}$ in
Assumption 1 is related to the dependence structure for the
lifetimes $T_{i}$, as will be further elaborated on below.

Next comes the main result, stating weak convergence on compact sets of the Nelson-Aalen estimator. 
\bth \label{thm.1}
Assume that $H$ is bounded and of finite variation, $H'(u-)(1-H(u-))^{-2}\in {\mathbf L}^{1}$ and that Assumption 1 holds. Then
\begin{eqnarray*}
      a_{n}^{-1}(\Lambda_{n}(t)-\Lambda(t))&\stackrel{\cal L}{\rightarrow}&
           \int_{0}^{t}\frac{1}{1-H(u-)}\,dw^{1}(u)+\int_{0}^{t}\frac{w(u)}{(1-H(u-))^{2}}\,dH^{1}(u),
\end{eqnarray*}
on ${\mathbf D}[0,\infty)$, as $n\rightarrow\infty$.
\eth
It is possible to obtain slightly more
general results in Theorem 1 (and Theorem 2 below), by allowing different rates
$a_{n}^{-1},b_{n}^{-1}$ for the two individual processes in Assumption 1; this would imply that only one of the two terms in Theorem
1 and 2 will turn up in the limit; we omit this generalization, note however the proof of Theorem 4 in the sequel. Such a generalization could be applicable to the (admittably strange) case when one of the processes is over-smoothed. 
\par
The Kaplan-Meier estimator of the survival function $S(t)$ is defined as the product-limit estimator
\begin{eqnarray}
      {S}_{n}(t)&=& \pi_{u\leq t}(1-d\Lambda_{n}(u)),\label{eq:prodlim} \\
      &=&\prod_{t_{i}\leq t} (1-\frac{\delta_{i}}{Y(t_{i})}), \nonumber
\end{eqnarray}
cf.$\!$ Gill and Johansen \cite{GillJohansen}. 
The limit distribution for the Kaplan-Meier estimator follows from limit process results for the Nelson-Aalen estimator coupled with the smoothness of the  map defining the Nelson-Aalen estimator, i.e. the product integral (cf. equation $(\ref{eq:prodlim})$)
\begin{eqnarray*}
      \phi:{\mathbf D}[0,\infty)\ni x(\cdot)&\mapsto& \pi_{u\leq \cdot}(1-dx(u))\in {\mathbf D}[0,\infty)
\end{eqnarray*}
is a compactly (or Hadamard) differentiable map, as shown in Proposition II.8.7 of Anderson et al.$\!$ \cite{ABGK}. That $\phi$ is compactly differentiable  means that for each $x\in {\mathbf D}[0,\infty)$ there is a linear functional $\phi'_{x}$ on ${\mathbf D}[0,\infty)$ such that
\begin{eqnarray*}
      \sup_{h\in K}||\frac{\phi(x+th)-\phi(x)}{t}-\phi'_{x}(h) ||&\rightarrow&0,
\end{eqnarray*}
as $t\rightarrow 0$, with $K$ an arbitrary compact set in ${\mathbf D}[0,\infty)$ and with $||\cdot||$ denoting the (supnorm over compact intervals) norm in ${\mathbf D}[0,\infty)$. Compact differentiability is precisely what is needed in statistical applications since it can be coupled with the process weak convergence results to yield weak convergence results for the studied estimator, via the functional delta method, cf.$\!$ Gill \cite{GillHadamard}. Compact differentiability of $\phi$ implies that if for some sequence $d_{n}\downarrow 0$ and stochastic process $v\in {\mathbf D}[0,\infty)$ we have
\begin{eqnarray*}
     d_{n}^{-1}(\Lambda_{n}-\Lambda)&\stackrel{\cal L}{\rightarrow}& v
\end{eqnarray*}
in $ {\mathbf D}[0,\infty)$, then the functional delta method implies that, if $S_{n}=\phi(\Lambda_{n})$,
\begin{eqnarray*}
     d_{n}^{-1}(S_{n}-S)&\stackrel{\cal L}{\rightarrow}& \phi'_{\Lambda}(v)
\end{eqnarray*}
in $ {\mathbf D}[0,\infty)$, as $n\rightarrow\infty$.
\bth \label{thm.2}
Assume that $H'(u-)(1-H(u-))^{-2}\in {\mathbf L}^{1}$ and that Assumption 1 holds. Then
\begin{eqnarray*}
        a_{n}^{-1}(S_{n}(t)-S(t))&\stackrel{\cal L}{\rightarrow}& [\int_{0}^{t}\frac{1}{1-H(u-)}\,dw^{1}(u)+\int_{0}^{t}\frac{w(u)}{(1-H(u-))^{2}}\,dH^{1}(u)]\cdot S(t),
\end{eqnarray*}
on ${\mathbf D}[0,\infty)$, as $n\rightarrow\infty$.
\eth
\section{Applications to different dependence structures}\label{Sec:3}
The Nelson-Aalen and Kaplan-Meier estimators' limit distributions are well known for independent data. In this section we recall these results, cf. van der Vaart \cite{vandervaart.2000}, as well as state generalization to different forms of dependence. 
\par
In the case of stationary and dependent data $\{t_{i}\}$ we can classify the results according to whether $\sum |\Cov(k)| < \infty$ or not, where $\Cov(k)=\Cov(t_{i},t_{i+k})$ is the covariance function. If $\sum |\Cov(k)|<\infty$ we say the sequence is weakly dependent and if not it is called strongly or long range dependent.
\par
Now assume that $\{t_{i}\}_{i\geq 1}$ is weakly dependent so that $\sum |\Cov(k)|<\infty$. Define the $\sigma$-algebras
\begin{eqnarray*}
        {\cal F}_{k}&=&\sigma\{t_{i}:i\leq k\},\\
        \bar{{\cal F}}_{k}&=&\sigma\{t_{i}:i\geq k\}.
\end{eqnarray*}
Then the sequence is called $\phi$-mixing if there is a function $\phi(n)\rightarrow 0$ as $n\rightarrow \infty$ and
\begin{eqnarray*}
        \sup_{A\in{\cal F}_{n}}|P(A|{\cal F}_{0})-P(A)|&\leq &\phi(n),
\end{eqnarray*}
almost surely.
\bth\label{thm.3}
Assume that $\{t_{i}\}$ is a stationary sequence with $E(t_{i})=0$ and $\sigma^{2}=\lim_{n\rightarrow\infty} n^{-1} \Var(t_{1}+\ldots +t_{n})>0$, and that either of
\begin{eqnarray*}
  (i)&&\{t_{i}\} \mbox{ are independent},\\
  (ii)&&\{t_{i}\} \mbox{ is $\phi$-mixing with }\sum \phi(k)^{1/2}<\infty,
\end{eqnarray*}
hold. Then if $H'(u-)(1-H(u-))^{-2}\in {\mathbf L}^{1}$ and Assumption 1 holds we obtain
\begin{eqnarray*}
        (n/\sigma^{2})^{1/2}(\Lambda_{n}(t)-\Lambda(t))&\stackrel{\cal L}{\rightarrow}& \int_{0}^{t}\frac{dB(H^{1}(u))}{1-H(u-)}+\int_{0}^{t}\frac{B(H(u))\,dH^{1}(u)}{(1-H(u-))^{2}} ,\\
        (n/\sigma^2)^{1/2}(S_{n}(t)-S(t))&\stackrel{\cal L}{\rightarrow}&[\int_{0}^{t}\frac{dB(H^{1}(u))}{1-H(u-)}+\int_{0}^{t}\frac{B(H(u))\,dH^{1}(u)}{(1-H(u-))^{2}}]\cdot  S(t),
\end{eqnarray*}
on ${\mathbf D}[0,\infty)$, as $n\rightarrow\infty$, with $B$ a standard Brownian bridge.
\eth
\bnote
Different combinations of mixing for $\{T_{i}\}$ and independence or mixing for $\{\tau_{i}\}$ are possible. We have stated the theorem in order to avoid technicalities concerning this. It is intuitively clear that mixing $\{T_{i}\}$ and independent or mixing $\{\tau_{i}\}$ result in mixing $\{t_{i}\}$, and so we refrain from developping this in more detail.
\enote
Note that the mixing part of Theorem 3 holds for other types of mixing condition also, what is necessary is a multivariate central limit theorem for dependent data, and we refrain from developing this in more detail. 
\par
In the long range dependent case, i.e. when $\sum |\Cov(k)|=\infty$, we are able to derive the statement of Assumption 1 only under rather strong assumptions. For this we review the limit distribution results for the empirical process for long range dependent data, of Dehling and Taqqu \cite{DehlingTaqqu}. Thus we assume that the terms in the generic process $\{\tau_{i}\}_{i\geq 1}$ can be written $\tau_{i}=g(\xi_{i})$ where $\{\xi_{i}\}_{i\geq 1}$ is some Gaussian stationary stochastic process with mean zero and covariance $\Cov_{\xi}(k)=\mbox{Cov}(\xi_{i},\xi_{i+k})=k^{-d}l_{0}(k)$ where $0<d<1$ is fixed, $l_{0}$ is a function slowly varying at infinity and $g$ is a function satisfying $E(g^{2}(\xi_{1}))$, cf.$\!$ Taqqu \cite{Taqqu1}. We call such a process $\{\tau_{i}\}_{i\geq 1}$ a subordinated Gaussian sequence (with parameter $d$). Next, for each fixed $t$, we can expand the terms in the empirical processes in a series in Hermite polynomials according to
\begin{eqnarray}
      1_{\{\tau_{1}\leq t\}}-F(t)&=&\sum_{k=r(t)}^{\infty}\frac{1}{k!}\eta_{k}(t)h_{k}(\xi_{1}), \label{eq:Hermiteexp}
\end{eqnarray}
with $h_{k}=$ the Hermite polynomial of order $k$, $\eta_{k}(t)=E[(1_{\{\tau_{1}\leq t\}}-F(t))h_{k}(\xi_{1})]$ the ${\mathbf L}^{2}$-projections of $1_{\{\tau_{1}\leq t\}}-F(t)$ on $h_{k}$, and with $r(t)$ the first nonzero coefficient in the expansion. We define the Hermite rank of the functions $\{ 1_{\{\tau_{1}\leq t\}}-F(t):t\in{\mathbf R}\}$ as $r=\inf_{t}r(t)$  and assume that $0<d<1/r$. Then also the sequence $\{ 1_{\{\tau_{i}\leq t\}}-F(t)\}$ is long range dependent and $\sigma_{n}^{2}=$Var$(\sum_{i=1}^{n}h_{r}(\xi_{i}))$ is asymptotically proportional to $n^{2-rd}l_{1}(n)$, where $l_{1}(n)=\frac{2}{r!(1-rd)(2-rd)}l_{0}(n)^{r}$. Furthermore, from Theorem 1.1 in Dehling and Taqqu \cite{DehlingTaqqu} it follows that
\begin{eqnarray*}
       \sigma_{n}^{-1}n(F_{n}(t)-F(t))&\stackrel{\cal L}{\rightarrow}& \frac{\eta_{r}(t)}{r!}z_{r}
\end{eqnarray*}
on ${\mathbf D}[0,\infty]$ (equipped with the supnorm metric), as $n\rightarrow\infty$, with $z_{r}$ a random variable.
\par
Note that in this case, the limit process $\eta_{r}(t)z_{r}/r!$ is stochastic only in that it has a stochastic scale function $z_{r}$; the process/function behaviour is described by the deterministic function $\eta_{r}(t)$. Note also that the limit is determined by the first nonzero coeffient in the Hermite series expansion $(\ref{eq:Hermiteexp})$. Further, the random variable $z_{1}$ is Gaussian while $z_{r}$ is not Gaussian for $r>1$, cf.$\!$ Taqqu \cite{Taqqu1}.
\par
In order to keep things simple we do not make seperate statements on the dependence of the two sequences $\{T_{i}\},\{\tau_{i}\}$, instead we assume a form of joint long range dependence in the next assumption. 
\bass
Assume that the sequences $\{t_{i}\}$ and $\{t_{i}\delta_{i}\}$ are both subordinated Gaussian long range dependent sequences, so that
\begin{eqnarray*} 
   t_{i}&=&g_{0}(\xi_{i}),\\
   t_{i}\delta_{i}&=&g_{1}(\xi_{i}),
\end{eqnarray*}  
where $\Var(g_{0}(\xi_{i})),\Var(g_{1}(\xi_{i}))<\infty$, and $\{\xi_{i}\}_{i\geq 1}$ is a Gaussian stationary stochastic process with mean zero and covariance $\Cov(\xi_{i},\xi_{i+k})=k^{-d}l_{0}(k)$ with $0<d<1$ is fixed, and $l_{0}$ a function slowly varying at infinity.
\eass
Note that the process $\{t_{i}\delta_{i}\}$ has a point mass in 0.
\bth \label{thm.4}
If Assumption 2 and $H'(u-)(1-H(u-))^{-2}\in {\mathbf L}^{1}$ hold, then
\begin{eqnarray*}
n^{rd/2}l_{1}(n)^{-1}(\Lambda_{n}(t)-\Lambda(t))&\stackrel{\cal L}{\rightarrow}& \zeta_{r}(t) \eta_{r} ,\\
       n^{rd/2}l_{1}(n)^{-1}(S_{n}(t)-S(t))&\stackrel{\cal L}{\rightarrow}& \zeta_{r}(t) \eta_{r}\cdot  S(t),
\end{eqnarray*}
on ${\mathbf D}[0,\infty)$, as $n\rightarrow\infty$. Here 
\begin{eqnarray*}
     \zeta_{r}(t)&=&\left\{ \begin{array}{ll} 
                           \int_{0}^{t}\frac{\eta_{r}(u)\,dH(u)}{(1-H(u-))^{2}}+ \int_{0}^{t}\frac{d\eta_{r}'(u)}{1-H(u-)}\mbox{ if }r_{0}\geq r_{1}    , \\ 
                            \int_{0}^{t}\frac{\eta_{r}(u)\,dH(u)}{(1-H(u-))^{2}} \mbox{ if } r_{0}<r_{1} , 
                                       \end{array} \right. \\ 
\end{eqnarray*}
and $r=\min(r_{0},r_{1})$, where $r_{0}$ is the Hermite rank of $1_{\{t_{i}\leq t,\delta_{i}=0\}}-H^{0}(t)$ and $r_{1}$ that of $1_{\{t_{i}\leq t,\delta_{i}=1\}}-H^{1}(t)$. Further
\begin{eqnarray*}
      \eta_{r}(t)&=&\left\{ \begin{array}{ll} 
                           \eta_{r,0}(t)+\eta_{r,1}(t)\mbox{ if }r_{0}\geq r_{1}    , \\ 
                             \eta_{r,0}(t)\mbox{ if } r_{0}<r_{1} , 
                                       \end{array} \right. \\ 
\end{eqnarray*}
with $\eta_{r,0}$ ($\eta_{r,1}$) the projection of $1_{\{t_{i}\leq t,\delta_{i}=0\}}-H^{0}(t)$ ($1_{\{t_{i}\leq t,\delta_{i}=1\}}-H^{1}(t)$) on the Hermite function.
\eth

Note that the limit processes in Theorem \ref{thm.4} are  degenerate, in that the function behaviour is described by deterministic functions with the random effect $\eta_r$ a factor. Note that $\eta_1$ is a Gausssian random variable, and $\eta_r$ is not for $r\geq 2$.

\section{Comments}
In this short note we have extended the limit distribution results for the Nelson-Aalen and Kaplan-Meier estimators for i.i.d. data to estimating the marginal distribution and corresponding cumulative hasard for stationary data, treating weakly dependent as well as long range dependent data.

We would like to stress the basic independence assumption: That the censoring times and life length times are independent, and thus that the censoring mechanism is independent of the life lengths for the events of interest.

\begin{appendix}
\section{Proofs}\label{App}
\prf (Lemma 1)
Let $\epsilon>0$ be given. The assumption of the lemma means that there are $K_{1},K_{2}$ compact subsets of ${\mathbf D}[0,\infty)$ such that $\inf_{n}P(w_{n}\in K_{1})>1-\epsilon/2$ and $\inf_{n}P(w_{n}^{1}\in K_{2})>1-\epsilon/2$. Define $K=K_{1}\times K_{2}$ and note that $K$ is compact in ${\mathbf D}[0,\infty)\times{\mathbf D}[0,\infty)$ by Tychonoff's theorem. Furthermore 
\begin{eqnarray*}
   \inf_{n} P((w_{n},w_{n}^{1})\in K)&=&\inf_{n}P(w_{n}\in K_{1},w_{n}^{1} \in K_{2})\\
        &\geq &\inf_{n}(P(w_{n}\in K_{1})+P(w_{n}^{1}\in K_{2})-1)\\
        &\geq & \inf_{n}P(w_{n}\in K_{1})+\inf_{n}P(w_{n}^{1}\in K_{2})-1\\
        &\geq & 1-\epsilon,
\end{eqnarray*}
and thus $(w_{n},w_{n}^{1})$ is tight in ${\mathbf D}[0,\infty)\times{\mathbf D}[0,\infty)$.
 
\prf (Theorem 1)
Note first that Assumption 1 implies that the individual processes $w_{n}$ and $w_{n}^{1}$ converge weakly to $w$ and $w^{1}$ respectiveley, as noted e.g.$\!$ in Billingsley \cite{Billingsley}.
\par
We have
\begin{eqnarray}
      a_{n}^{-1}(\Lambda_{n}(t)-\Lambda(t))&=&a_{n}^{-1}(\int_{0}^{t}\frac{dH_{n}^{1}(u)}{1-H_{n}(u-)}-\int_{0}^{t}\frac{dH^{1}(u)}{1-H(u-)}).\label{eq:bevis1}
\end{eqnarray}
Keep $u$ fixed, and make a series expansion of $(1-H_{n}(u-))^{-1}$ around $H(u-)$ to obtain
\begin{eqnarray}
     \frac{1}{1-H_{n}(u-)}&=&\frac{1}{1-H(u-)}+\frac{1}{(1-H(u-))^{2}}(H_{n}(u-)-H(u-))(1+o(1))\label{eq:seriesexp}
\end{eqnarray}
as $H_{n}(u-)\rightarrow H(u-)$. But from Assumption 1 it follows that for all compact sets $K\subset [0,\infty)$ (and so in particular for the compact set $[0,t]$), we have $\sup_{u\in K}|H_{n}(u)-H(u)|=o_{P}(1)$ as $n\rightarrow\infty$, which implies that $(\ref{eq:seriesexp})$ actually holds uniformly over all $u\in[0,t]$, if we replace the $1+o(1)$ term with $1+o_{P}(1)$. This implies
\begin{eqnarray}
     a_{n}^{-1}(\Lambda_{n}(t)-\Lambda(t))&=&\int_{0}^{t}\frac{1}{(1-H(u-))}\,a_{n}^{-1}d(H_{n}^{1}-H^{1})(u)\nonumber\\
     &&+\int_{0}^{t}\frac{a_{n}^{-1}(H_{n}(u-)-H(u-))}{(1-H(u-))^{2}}\,dH_{n}^{1}(u)(1+o_{P}(1))\nonumber\\
     &=:&I_{1}(t)+I_{2}(t).\label{eq:bevis2}
\end{eqnarray}
Note that here the $1+o_{P}(1)$ holds uniformly for all $t$ in compact intervals $[0,s]$: since $1+o_{P}(1)$ holds uniformly for all $u$ in compact sets and since $\sup_{[0,s]}|\int_{0}^{t} f(u)\,du|\leq s\sup_{[0,s]}|f(s)|$.
\par
Clearly, to show weak convergence of the term $I_{1}$, it is enough to show that the (marginal) map  
\begin{eqnarray*}
    U:{\mathbf D}([0,\infty),{\mathbf R})\times {\mathbf D}([0,\infty),{\mathbf R})\ni (v,w)&\mapsto& \int_{0}^{t}\frac{1}{1-H(u-)}\,dv(u)\in {\mathbf C}[0,\infty)
\end{eqnarray*}
is continuous. By partial integration, allowed if $v$ is of finite variation, we have 
\begin{eqnarray*}
      U(v,w)&=&v(t)\frac{1}{1-H(t-)}-\int_{0}^{t}\frac{H'(u-)}{(1-H(u-))^{2}}v(u)\,du,
\end{eqnarray*}
 assuming (without loss of generality) that $v(0)=0$. Thus if $H'(u-)(1-H(u-))^{-2}\in {\mathbf L}^{1}$ clearly $v_{n}\rightarrow v$ (in supnorm on compact intervals) will imply that the second term of $U(v_{n})$ converges to the second term of $U(v)$ (in supnorm on compact intervals). For the first term of $U(v_{n})$ the convergence to the first term of $U(v)$ (uniformly on compact intervals) holds trivially, if $(1-H(t-)^{-1}$ is bounded on compact intervals. Thus $U$ is a continuous map.
\par
For the term $I_{2}$ we have
\begin{eqnarray*}
  I_{2}(t)&=&\int_{0}^{t}\frac{a_{n}^{-1}(H_{n}(u-)-H(u-))}{(1-H(u-))^{2}}\,dH^{1}(u)(1+o_{P}(1))\\
  &&+\int_{0}^{t}\frac{H_{n}(u-)-H(u-)}{(1-H(u-))^{2}}\,a_{n}^{-1}d(H_{n}^{1}-H^{1})(u)(1+o_{P}(1)).
\end{eqnarray*}
To show weak convergence of the first term of $I_{2}$, it is enough to show continuity of the (marginal) map
\begin{eqnarray*}
 V:{\mathbf D}([0,\infty),{\mathbf R})\times {\mathbf D}([0,\infty),{\mathbf R})\ni (u,v)&\mapsto& \int_{0}^{t}\frac{v}{(1-H(u-))^2}\,dH^{1}(u)\in {\mathbf C}[0,\infty),
\end{eqnarray*}
which is clear, since $H^{1'}(u)(1-H(u-))^{-2}\in{\mathbf L}^{1}$ if $H'(u)(1-H(u-))^{-2}\in{\mathbf L}^{1}$. The second term of $I_{2}$ can be bounded, uniformly in compact intervals $[0,k]$, by
\begin{eqnarray*}
     \sup_{t\in[0,k]} |H_{n}(t)-H(t)|\sup_{t\in[0,k]}|\int_{0}^{t}\frac{1}{(1-H(u-))^{2}}|a_{n}^{-1}d(H_{n}^{1}-H^{1})(u)|(1+o_{P}(1))
     &=&O_{P}(a_{n}),
\end{eqnarray*}
as $n\rightarrow\infty$, implying that this term is negligible.
\par
Thus we finally obtain
\begin{eqnarray*}
  a_{n}^{-1}(\Lambda_{n}(t)-\Lambda(t))&=&U(a_{n}^{-1}(H_{n}^{1}-H^{1}))+ V(a_{n}^{-1}(H_{n}-H))+O_{P}(a_{n})
\end{eqnarray*}
(as $n\rightarrow\infty$, uniformly for $t$ in compact intervals), which by the continuous mapping theorem, Slutsky's theorem and since addition of the marginals ${\mathbf D}([0,\infty),{\mathbf R})\times {\mathbf D}([0,\infty),{\mathbf R})\mapsto{\mathbf D}[0,\infty) $ is continuous, implies weak convergence
\begin{eqnarray*}
      a_{n}^{-1}(\Lambda_{n}-\Lambda)&\stackrel{\cal L}{\rightarrow}&U(w^{1})+V(w^{2}),
\end{eqnarray*}
on ${\mathbf D}[0,\infty)$, which implies the statement of the theorem. 
 
\prf (Theorem 2)
The only thing that remains is to derive the functional derivative of the product integral $\phi$. This is done in Proposition II.8.7 in Anderson et al.$\!$ \cite{ABGK} (note that we have scalar $S$ and $\Lambda$ so that multiplication is commutative in our case) from which it follows that
\begin{eqnarray*}
         \phi'_{\Lambda}(v)(t)&=&\prod_{0<u\leq t}(1-d\Lambda(u))\cdot v(t)=S(t)v(t).
\end{eqnarray*}
Thus an application of Theorem 1 and the functional delta method, cf.$\!$ Gill \cite{GillHadamard}, ends the proof.
 
\prf (Theorem 3)
We will prove the statement of Assumption 1 with $a_{n}=n^{1/2}$, $w^{1}(t)=(\sigma^2)^{1/2}B(H(t))$ and $w^{2}(t)=(\sigma^2)^{1/2}B(H^{1}(t))$, and with $B$ a standard Brownian bridge. For this we will show convergence of the finite dimensional distribution of $\{w_{n}\}_{n\geq 1}$ and tightness. By the Cram\`er-Wold device, convergence of the finite dimensional distributions of the bivariate process $(w_{n},w_{n}^{1})$ is equivalent to weak convergence of the two-dimensional vector
\begin{eqnarray*}
  &&\gamma_{1}(w_{n}(t_{1}),w_{n}^{1}(t_{1}))+\ldots +\gamma_{k} (w_{n}(t_{k}),w_{n}^{1}(t_{k}))\\
  &=&(\gamma_{1}w_{n}(t_{1})+\ldots \gamma_{k}w_{n}(t_{k})   ,\gamma_{1}w_{n}^{1}(t_{1})+\ldots \gamma_{k}w_{n}^{1}(t_{k})    )
\end{eqnarray*}
which is equivalent to weak convergence of the random variable
\begin{eqnarray*}
    \lambda_{1}(\gamma_{1}w_{n}(t_{1})+\ldots \gamma_{k}w_{n}(t_{k})  )+\lambda_{2}( \gamma_{1}w_{n}^{1}(t_{1})+\ldots \gamma_{k}w_{n}^{1}(t_{k})   )   
\end{eqnarray*}
which implies that it is enough to prove weak convergence of the univariate stochastic process $w_{n}^{\lambda}=w_{n}+\lambda w_{n}^{1}$ for every fixed $\lambda$, for $i.i.d.$ and $\phi-$mixing data (cf. Philipp \cite{Philipp}) and since addition of the marginals in ${\mathbf D}([0,\infty),{\mathbf R})\times {\mathbf D}([0,\infty),{\mathbf R})$ is a continuous map.
\par
Since $w_{n}^{\lambda}$ is a normalized sum of mean zero simple functions, by the multivariate central limit theorem it follows that the finite dimensional distributions of $w_{n}^{\lambda}$ converge to those of $(\sigma^2)^{1/2}[B(H(t))+\lambda B(H^{1}(t))]$, where $B$ is a standard Brownian bridge.
\par
Tightness follows from Lemma 1, since $\{w_{n}\}_{n\geq 1}$ and $\{w_{n}^{1}\}_{n\geq 2}$ are tight, and the theorem is proven.
 
\prf (Theorem 4)
We will prove the statement of Assumption 1, with $a_{n}=\sigma_{n}n^{-1}\sim n^{-rd/2}l_{1}(n)$. As in the proof of Theorem 3, we need to prove that the finite dimensional distributions of $w_{n}^{\lambda}=w_{n}+\lambda w_{n}^{1}$ converge and that $\{(w_{n},w_{n}^{1})\}_{n\geq 1}$ is tight. Convergence of the finite dimensional distributions is obtained by proving that $w_{n}^{\lambda}$ converges weakly for each $\lambda$. Since 
\begin{eqnarray*}
          w_{n}^{\lambda}(t)&=&\frac{1}{n}\sum_{i=1}^{n}(1_{\{t_{i}\leq t\}}+\lambda\delta_{i}1_{\{t_i\leq t\}}-H(t)-\lambda H^{1}(t)),
\end{eqnarray*}
we can make a series expansion of the terms of $w_{n}^{\lambda}(t)$ as 
\begin{eqnarray*}
        1_{\{t_{i}\leq t\}}+\lambda\delta_{i}1_{\{t_i\leq t\}}-H(t)-\lambda H^{1}(t)&=&\sum_{k=r(t)}^{\infty}\frac{\eta_{k}^{\lambda}(t)}{k!}h_{k}.
\end{eqnarray*}
Clearly, since projection is a linear operation, we have that $\eta_{k}^{\lambda}(t)=\eta_{k,0}(t)+\lambda \eta_{k,1}(t)$, where $\eta_{k,0}(t),\eta_{k,1}(t)$ are defined in the statement of the theorem. 
\par
Let $r^{\lambda},r_{0},r_{1}$ be the Hermite ranks corresponding to $\eta_{k}^{\lambda},\eta_{k,0},\eta_{k,1}$ respectively. Assume first that $r_{0}\geq r_{1}$; then $H_{n}$ and $H_{n}^{1}$ will have the same rate of convergence. From Theorem 1.1 in Dehling and Taqqu \cite{DehlingTaqqu}, it follows that $a_{n}^{-1}w_{n}^{\lambda}$ converges weakly to $(\eta_{r,0}(t)+\lambda\eta_{r,1}(t))z_{r}/r!$.
\par
To prove tightness, clearly $\{w_{n}\}_{n\geq 1}$ and $\{w_{n}^{1}\}_{n\geq 2}$ are separately tight, since they have weak limits by Theorem 1.1 in Dehling and Taqqu \cite{DehlingTaqqu}, and thus from Lemma 1 $\{(w_{n},w_{n}^{1})\}_{n\geq 1}$ is tight. We have proven Assumption 1, i.e. that
\begin{eqnarray*}
        n^{rd/2}l_{1}(n)^{-1}(w_{n}(t),w_{n}^{1}(t))&\stackrel{\cal L}{\rightarrow}& \frac{1}{r!}(\eta_{r,0}(t),\eta_{r,1}(t))z_{r},
\end{eqnarray*}
on ${\mathbf D}([0,\infty),{\mathbf R})\times {\mathbf D}([0,\infty),{\mathbf R})$, as $n\rightarrow\infty$. 
\par
If instead the Hermite ranks satsify $r_{0}<r_{1}$ then instead $H_{n}$ and $H_{n}^{0}$ will have the same rate of convergence, and then $H_{n}^{1}$ will be negligible. In that case
\begin{eqnarray*}
        n^{rd/2}l_{1}(n)^{-1}(w_{n}(t),w_{n}^{1}(t))&\stackrel{\cal L}{\rightarrow}& \frac{1}{r!}(\eta_{r,0}(t),0)z_{r},
\end{eqnarray*}
on ${\mathbf D}([0,\infty),{\mathbf R}\times {\mathbf D}([0,\infty),{\mathbf R}))$, as $n\rightarrow\infty$. An application of Theorem 1 and 2 ends the proof.

\end{appendix}
\par

\section*{References}

\bibliography{referencefile}

\end{document}